\documentclass[notitlepage, 11pt]{article}

\usepackage{color}
\usepackage{latexsym}
\usepackage{amsmath, amssymb, amsfonts}
\usepackage[toc]{appendix}
\numberwithin{equation}{section}

\newtheorem{lemma}{Lemma}
\newtheorem{theorem}{Theorem}

\newtheorem{proposition}{Proposition}

\newtheorem{assumption}{Assumption}

\newtheorem{remark}{Remark}

\newcommand{\beginsec}{
\setcounter{lemma}{0}
\setcounter{theorem}{0}
\setcounter{corollary}{0}
\setcounter{definition}{0}
\setcounter{example}{0}
\setcounter{proposition}{0}
\setcounter{condition}{0}
\setcounter{assumption}{0}
\setcounter{conjecture}{0}
\setcounter{problem}{0}
\setcounter{remark}{0}
}

\newcommand{\noi}{\noindent}

\newcommand{\E}{\mathbb{E}}
\newcommand{\R}{\mathbb{R}}

\newcommand{\N}{\mathbb{N}}

\newcommand{\Ir}{\mathbb{I}}

\newcommand{\eps}{\varepsilon}

\newcommand{\ph}{\varphi}

\newcommand{\PP}{{\mathbb P}}

\newcommand{\calA}{{\cal A}}

\newcommand{\calC}{{\cal C}}

\newcommand{\calE}{{\cal E}}
\newcommand{\calF}{{\cal F}}

\newcommand{\calH}{{\cal H}}

\newcommand{\bpin}{\boldsymbol{\pi^n}}
\newcommand{\bpins}{\boldsymbol{\pi^{*,n}}}
\newcommand{\bpis}{\boldsymbol{\pi^*}}

\newcommand{\bPi}{\boldsymbol{\Pi}}
\newcommand{\bpi}{\boldsymbol{\pi}}
\newcommand{\bkappan}{\boldsymbol{\kappa^n}}

\newcommand{\skp}{\vspace{\baselineskip}}

\newcommand{\iy}{\infty}

\newcommand{\ds}{\displaystyle}
\newcommand{\A}{{\cal A}}

\newcommand{\Pu}{\rho}

\title{Risk Sensitive Control of the Lifetime Ruin Problem}
\author{Erhan Bayraktar\thanks{web: www-personal.umich.edu/$\sim$ erhan/, email: erhan@umich.edu} \hspace{4em} Asaf Cohen\thanks{web: https://sites.google.com/site/asafcohentau/, email: asafc@umich.edu}\\ \\ \\
Department of Mathematics\\
University of Michigan\\
Ann, Arbor, 48109, USA}

\date{\today}

\begin{document}

\maketitle

\begin{abstract}

\skp
We study a risk sensitive control version of the lifetime ruin probability problem. We consider a sequence of investments problems in Black-Scholes market that includes a risky asset and a riskless asset. We present a differential game that governs the limit behavior. We solve it explicitly and use it in order to find an asymptotically optimal policy.

\skp

\noi{\bf Keywords:} Probability of lifetime ruin, optimal investment, risk sensitive control, large deviations, differential games. \,\,
\end{abstract}

\section{Introduction}
\beginsec
%
%

The problem of how an individual should invest her wealth in a risky financial market in order to minimize the probability of outliving her wealth, also known as the probability of lifetime ruin was extensively analyzed, see e.g.~\cite{Milevsky2000}, \cite{Young2004}, \cite{Bayraktar2007}, \cite{Bayraktar2007b}, \cite{Bayraktar2010}, \cite{Bayraktar2011b}, \cite{Yener2014}, and \cite{bayraktar2015minimizing}. These works fall naturally within the area of optimally controlling wealth to reach a goal. Research on this topic goes back to the seminal work of \cite{Dubins1965} and continued with the work of \cite{Pestien1985}, \cite{Orey1987}, \cite{Sudderth1989}, \cite{Kulldorff1993}, \cite{Karatzas1997}, \cite{Browne1997}, \cite{Browne1999a}, and \cite{Browne1999b}.

In the standard Black-Scholes market that includes a risky asset and a riskless asset, the case of interest is when the investor consumes more than the potential profit that follows by investing the entire wealth in the riskless asset, that is $c(x)>rx$, in which $c(\cdot)$ is the consumption function, $r$ is the constant riskless rate, and $x$ is the current wealth. The other case is trivial, of-course, since by investing the entire wealth in the riskless asset the wealth cannot decrease and ruin is avoided. In case that $c(x)-rx\approx 0^+$ then the investor who wishes to minimize the probability of lifetime ruin should invest almost all of her wealth in the riskless asset. The probability of ruin would be small, yet positive. With the understanding that lifetime ruin is a rare and dramatic event and that one should also avoid living close to the ruin level, we study this case, by using a risk sensitive control framework. The risk sensitive control criteria, penalizes such events heavily, and therefore, provides a natural way to address these considerations.

We study the risk sensitive control via large deviations techniques. In \cite{Pham2007}, Pham provides some applications and methods of large deviations in finance and insurance. Among the studied models, he considers ruin problems when the initial reserve is large and therefore, the probability of ruin is small. We, on the other hand, study a lifetime ruin problem, which is a different problem, and via risk sensitive control with small noise, as described below, which yields a different analysis.

In order to rigorously treat the mentioned case that $c(x)-rx\approx 0^+$ we consider a sequence of models, indexed by $n\in\N$, that differ from each other only in the consumption function in a way that $c^n(x)-rx=O(1/n)$, where $n$ is a large parameter. By using an appropriate time scaling we get a risk sensitive control with small noise as follows: The scaled wealth process under the consumption function $c^n$ satisfies
\begin{align}\notag
d\tilde W^n(t) &= b(\tilde W^n(t),\bpin(t))dt+\frac{1}{\sqrt{n}}\sigma(\tilde W^n(t),\bpin(t))dB(t),\quad t\ge 0,\\\notag
\tilde W^n(0)&= x
\end{align}
for some proper $b$ and $\sigma$, where $\bpin$ is the investment policy, and $B$ is a standard Brownian motion. The goal is to choose $\bpin$ that minimizes
\begin{align}\notag
\frac{1}{n}\ln \E \left[  e^{n\left( \int_0^{\tau^n_a\wedge\tau^n_d} l(\tilde W^n(s))ds + \Pu 1_{\{\tau^n_a\leq
\tau^n_d\}}\right)}\right],
\end{align}
where $\tau^n_d$ is the time of death, $\tau^n_a$ is the time of reaching the ruin level $a$, $\Pu$ is a penalty for lifetime ruin, and $l$ is a nonnegative non-increasing Lipschitz function that penalizes low wealth. We present a differential game that governs the limiting behavior. We solve it explicitly and use it in order to find an asymptotically optimal policy.

Risk sensitive control for controlled stochastic differential equations with small noise have been studied for example in \cite{Fleming1995}, \cite{Fleming2006}, and \cite{Dupuis1989}. For a survey about the topic the reader is referred to \cite{Fleming2006b}. There are several approaches towards this problem. In \cite{Fleming2006}, Fleming and Soner used differential equations tools and show that the sequence of the appropriate prelimit Hamiltonians converges to the Hamiltonian that is associated with the differential game.
Among other requirements, it is assumed that the terminal cost is continuous and that the terminal time is fixed. In our case, the indicator takes the role of the terminal cost, which besides of being not continuous, in this case it also depends on the history of the wealth process. Also, we consider a random terminal time that is independent of the wealth process. Moreover, partial differential equations techniques does not provide asymptotically optimal policies, while we do.

In \cite{Dupuis1989}, Dupuis and Kushner approached a risk sensitive control problem of minimizing escape time probabilities by techniques taken from the theory of large deviation. Some of their requirements are that the drift and the diffusion coefficients, $b$ and $\sigma$ respectively, are bounded and the latter is also non-degenerate and does not depend on the control. These requirements are essential for the proofs. Also, they use a fixed terminal time. In our model, besides that the terminal time is random, the drift and the diffusion coefficients are assumed to be Lipschitz, but only the diffusion coefficient, $\sigma$, is assumed to be bounded. We allow $\sigma$ to be zero and to depend on the control. In fact, under the asymptotically optimal policy that we suggest the diffusion coefficient can be degenerate.

Recently, in \cite{Atar2014} and \cite{Atar2015b} the authors considered a queuing network problem under the moderate deviation heavy traffic regime. By using a variant of the proof of Varadhan's lemma and some properties of the differential game, an asymptotic optimality in the queueing systems is shown. In these papers, the controlled stochastic processes are not diffusion, but they are relatively close in distribution to a controlled diffusion with small noise. Therefore, the analysis requires some additional tools, and mainly the Skorohod mapping. While the structure of the queueing network in the prelimit raises some difficulties, the approximated diffusion is relatively simple and consists of Brownian motion (reflected Brownian motion, in the second paper) with drift. 
Although our proof considers some measure change arguments and is inspired by the proof of Varadhan's lemma,  in contrast to \cite{Atar2014} and \cite{Atar2015b}, we need to work with a controlled diffusion process. 

Regarding the random terminal time, the cost function can be referred as a discounted version of the risk sensitive cost. The only model from the above that considered a similar discounted structure is \cite{Atar2015b}. However, unlike the mentioned paper, we consider a scaled discount factor. The differential game associated with \cite{Atar2015b} appears in \cite{Atar2015a} and like in our case, the optimal solution of the game is time-homogeneous. Motivated by this property we analyze discounted risk sensitive control with small noise diffusions further in a future paper.

Let us summarize the contribution of this paper:
\begin{itemize} \itemsep0em
\item We propose a risk-sensitive cost for a lifetime ruin problem, which can be expressed as a discounted risk sensitive cost. We present a differential game that governs the limiting behavior.
\item We solve the differential game explicitly, including finding an optimal policy for the minimizer that leads to an asymptotically optimal policy in the prelimit stochastic model.
\item Our assumptions over the diffusion process are weaker than what usually appears in the literature, and yet 
we manage to find an asymptotically optimal control.
\end{itemize}

The organization of the paper is as follows. In Section \ref{sec2} we
describe the model, introduce the differential game, and state the main results.
In Section \ref{sec_3} we analyze the differential game,
present an Hamilton--Jacobi--Bellman (HJB) equation, characterize the
differential game's value function as its unique solution,
and we provide an explicit expression for the value function. Then we
present an explicit optimal control for the minimizer, and a simple control for the maximizer that achieves the value function.
In Section \ref{sec_4} we prove the main result by showing that in the
limit the differential game describes the stochastic model. 

We close this section by introducing some frequently used notation.

\noindent \textbf{Notation.}  
We denote $[0, \iy)$ by $\R_+$. For $f:[0,t]\to\R$ let $|f|_t:=\sup_{0\le s\le t}|f(s)|$.
For any interval $I$ denote by $\calA\calC (I)$ and $\calC (I)$ the
spaces of absolutely continuous functions (resp.,
continuous functions) mapping $I\to\R$.
Write $\calA\calC_0 (I)$ and $\calC_0 (I)$ for the subsets of the corresponding
function spaces, of functions that start at zero.

\section{Model and results}\label{sec2}
\beginsec

\subsection{The stochastic model}\label{sec2a}

We consider a sequence of stochastic models, indexed by $n\in\N$ of an investor who trades continuously in a Black-Scholes type financial market with no transaction costs. We allow borrowing and short-selling. The price of the riskless asset follows 
\begin{align}\notag
dV(t)=rV(t)dt,
\end{align}
where $r\ge 0$ is the constant interest rate. The risky asset follows a geometric Brownian motion:
\begin{align}\notag
dS(t) &= S(t)\left[\mu dt + \sigma dB(t)\right],
\end{align}
where $\mu>r$ and $\sigma>0$ are constants and $(B(t))_{t\ge 0}$ is a standard Brownian motion. 
For reasons that will be clear onwards we define a sequence of consumption functions, indexed by $n$. For any given $n\in\N$ we assume that consumption function takes the form: $c^n(\cdot)=r\cdot+\frac{1}{n}e(\cdot)$, for some function $e:[a,\infty)\to\R$. We assume that $e(\cdot)$ is a Lipshcitz function and that there is a positive constant $M_0$ such that $e(\cdot)\leq M_0$. 
For every $n\in\N$ and at any given time $t\ge 0$ let $\bkappan(t)$ be the amount of money that is invested in the risky asset. Then the wealth process satisfies
\begin{align}\notag
dW^n(t) &= \left(rW^n(t)-c^n(W^n(t))+(\mu-r)\bkappan(t)\right)dt+\sigma\bkappan(t)dB(t),\quad t\ge 0,\\\notag
W^n(0)&= x.
\end{align}
Now, by using time scaling and by referring to $\tilde W^n(\cdot)=W^n(n\cdot)$ we get that
\begin{align}\label{eq_203}
d\tilde W^n(t) &= \left(-e(\tilde W^n(t))+(\mu-r)\bpin(t)\right)dt+\frac{1}{\sqrt{n}}\sigma\bpin(t)dB(t),\quad t\ge 0,\\\notag
\tilde W^n(0)&= x,
\end{align}
where $\bpin(\cdot)=n\bkappan(n\cdot)$. In what follows we will denote by $\{\calF_t\}_{0\le t\le 2T}$, the usual augmentation of the natural filtration generated by the Brownian motion in \eqref{eq_203}.
From now onwards, we refer to $\bpin$ as the control. We denote by
$\bPi=\bPi_{M_1}$
the collection of all progressively measurable processes $(\bpi(t))_{t\ge 0}$ such that $|\bpi(\cdot)|\le M_1$, which we refer to as {\it admissible} policies, where $M_1$ is a positive constant. We take $\bpin\in\bPi$.
By the assumptions on $e(\cdot)$ and $\bpin$ it follows that for every $x>0$, the above admits a unique solution.
For every $n\in\N$, denote by $\tau^n_a$ the first time that $\tilde W^n$ reaches $a\in(0,x)$, which we will refer to as the {\it ruin level}. The investor would like to avoid ruin during her lifetime and also to avoid long living close to the ruin level. Also, let $\tau^n_d$ be the investor's random time of death. Due to the time scaling, we assume that $\tau^n_d$ is exponentially distributed with parameter $\lambda n$. The goal of the investor is to minimize the following risk sensitive control cost:

\begin{align}\notag
J^n(x,\bpin):&=\frac{1}{n}\ln \E \left[  e^{n\left( \int_0^{\tau^n_a\wedge\tau^n_d} l(\tilde W^n(s))ds + \Pu 1_{\{\tau^n_a\leq \tau^n_d\}}\right)}\right]\\\notag
&=\frac{1}{n}\ln \E \left[\int_0^\infty e^{-\lambda n t }  e^{n\left( \int_0^{\tau^n_a\wedge t} l(\tilde W^n(s))ds + \Pu 1_{\{\tau^n_a \leq t\}}\right)}dt\right] + \frac{1}{n}\ln(\lambda n),
\end{align}
where $\Pu>0$ stands for the punishment cost for ruining and $l:[a,\iy)\to [0,\lambda)$ is a non-increasing function. The function $l$ represents a punishment for the investor when her wealth is close to the ruin level $a$.  Obviously, we would like to give higher punishment when the wealth is closer to $a$. 
Moreover, since we would like that, given $n$, the function $J^n$ would be decreasing with respect to (w.r.t.) the wealth, we require that $l(\cdot)<\lambda$. Otherwise, $J^n$ would be increasing around $a$. This case represents a situation when the punishment of living close to the ruin level dominates the punishment from being ruined. Notice also that the last term on the right-hand side (r.h.s.) of the above goes to zero as $n$ goes to infinity. For this reason we will ignore it in the analysis in Section \ref{sec_4}. We summarize the assumptions mentioned  above:
\begin{assumption}\label{assumption1}
Set the following constants $\mu>r>0$, $\sigma,M_0,a,\lambda>0$. Also, let $e:[a,\iy)\to [-\iy,M_0]$ be Lipschitz and $l:[0,\iy)\to[0,\lambda)$ be a Lipschitz and non-increasing function. 
\end{assumption}
The assumption is at force throughout the paper. 

We study the problem when $n\to\iy$. 
As mentioned in the introduction, both the prelimit stochastic model and the limit suffer from several complexities in the analysis. First, the indicator part of the cost function complicates the analysis because it depends on the history of the process and if we look at it as a terminal cost, then it is not continuous w.r.t.~the terminal wealth. Second, we study a discounted version of the risk sensitive cost. To the best of our knowledge such formulation studied before only in \cite{Atar2015b} and also in a queueing system framework and with a discount that is free of $n$. Third, the diffusion coefficient is not necessarily bounded away from zero and it depends on the control. In fact as is shown in Section \ref{sec_2c}, the asymptotically optimal policy may become zero and therefore, so does the volatility coefficient. Therefore, the Hamiltonian method of \cite[Chapter XI.7]{Fleming2006}, or change of measure method of \cite{Dupuis1989} do not work here.  We find an asymptotically optimal policy for the problem by studying a differential game. We show that as $n\to\iy$ the optimal risk sensitive cost function  converges to the value of the game, and that an asymptotically optimal policy can be deduced from the minimizer's optimal control in the game.

%

\subsection{Differential game setting}\label{sec2b}

In this section, inspired by \cite[Chapter XI.7]{Fleming2006} and \cite[Theorem 5.6.7]{Dembo1998} we describe a differential game associated with the optimal risk sensitive control problem. We denote by $\Pi=\Pi_{M_1}$ the set of all Lipschitz functions $\pi:[a,\iy]\to[-M_1,M_1]$. Given $\pi\in\Pi$ and $\psi\in\calA\calC_0[0,\iy)$, the \textit{state process} associated with the initial condition $x$ and the data $\psi$ and $\pi$ is given by
\begin{align}\label{eq_250}
\dot\ph(t)&=-e(\ph(t))+(\mu-r)\pi(\ph(t)) + \sigma\pi(\ph(t))\dot\psi(t),\;\;t\geq 0,\\\notag
\ph(0)&=x.
\end{align}
One can easily verify that the state process is well-defined, see \cite[Theorem 19.12]{Poznyak2009}.
Note the analogy between the above and \eqref{eq_203}. The game payoff is
\begin{align}\notag
\sup_{T\in[0,\iy)}\Big\{\int_0^{T\wedge\tau}[-\lambda +l(\ph(t))]dt-\Ir(T\wedge\tau,\psi)+ \Pu 1_{\{\tau\leq T\}}\Big\},
\end{align}
where $\tau$ is the first time that the state process hits the ruin level $a$ and for every $t>0$, $\Ir(t,\cdot)$ is a function mapping
$\calC[0,t]$ to $\R_+\cup\{+\iy\}$ defined as
\begin{equation}\notag
\Ir(t,\psi):=\begin{cases}
                \ds
                \frac{1}{2}\int_0^t\dot\psi^2(s)ds &\mbox{if}\ \psi\in\calA\calC_0[0,t],
\\ \\
                +\infty & \mbox{otherwise}.
              \end{cases}
\end{equation}
The function $\Ir$ is the rate function\footnote{Although we are not using explicitly large deviation arguments in the paper, for intuition reasons we still choose to define the cost by using the rate function instead of simply using only $\tfrac{1}{2}\int_0^t\dot\psi^2(s)ds$.} of the Brownian motion $(\frac{1}{\sqrt{n}}B(t))_t$, as $n\to\iy$, see \cite[Theorem 5.2.3]{Dembo1998}. The ``$\sup_{T\in[0,\iy)}$'' is the differential game analogue of the control problem's discount factor
, $\lambda n$.
The payoff is maximized over $\psi$ and minimized over $\pi$. By the definition of the function $\Ir$ we may restrict the maximizer only to $\psi\in\calA\calC_0[0,\iy)$. 

The control $\pi\in\Pi$ is taken to be a feedback control and $\psi\in\calA\calC_0[0,\iy)$ and $T\in\R_+$ are open-loop controls. We call $\psi$ the \emph{path part of the control} and the $T$ a \emph{termination time  part of the control}.
Given $x\in[a,\infty)$, $\pi\in \Pi$, $\psi\in \calA\calC_0[0,\iy)$, and $T\in\R_+$, we define the cost until time $T$ by
\begin{align}\label{eq_253}
C(x,\pi,\psi,T):=\int_0^{T\wedge\tau}[-\lambda +l(\ph(t))-\frac{1}{2}\dot\psi^2(t)]dt +\Pu 1_{\{\tau\leq T\}}.
\end{align}
The value of the game is defined by
\begin{equation}\label{eq:defn-U}
U(x):=\inf_{\pi\in\Pi}\;\sup_{\psi\in \calA\calC_0[0,\iy),T\in\R_+}\;C(x,\pi, \psi,T).
\end{equation}
In the remark below we show that the maximizer can be restricted to a smaller set of controls without any loss. This property serves us in the sequel.

\begin{remark}\label{rem_21}
(1) Since $l(\cdot)<\lambda$ it follows that for every $\psi\in\calA\calC_0[0,\iy)$ and every $T\in\R_+$ one has $\int_0^{T\wedge\tau}[-\lambda +l(\ph(t))-\frac{1}{2}\dot\psi^2(t)]dt\leq 0$. Therefore, without any loss for the maximizer she can be restricted to $\psi$'s under which $\tau<\iy$ and $T \in \{0,\tau\}$.

\skp\noi
(2) Notice moreover that the maximizer can also be restricted to $\psi$'s for which the state process satisfies $\ph(t)<\ph(0)=:x$ for every $t>0$ and by (1) above also $\tau<\iy$. Indeed, since the integrand on the r.h.s.~of \eqref{eq_253} is negative then the only way that $U$ is positive is in case that $\tau<\iy$. Let $\psi=\psi_\pi$ be such that $\tau<\iy$. Denote by $\tau_x$ the last time before time $\tau$ that $\ph(t)=x$. Then,
\begin{align}\notag
C(x,\pi,\psi,\tau)&=\int_0^{\tau}[-\lambda +l(\ph(t))-\frac{1}{2}\dot\psi^2(t)]dt +\Pu \\\notag
&<\int_{\tau_x}^{\tau}[-\lambda +l(\ph(t))-\frac{1}{2}\dot\psi^2(t)]dt +\Pu \\\notag
&=\int_{0}^{\tau-\tau_x}[-\lambda +l(\ph_x(t))-\frac{1}{2}(\dot\psi_x)^2(t)]dt +\Pu \\\notag
&=C(x,\pi,\psi_x,\tau-\tau_x),
\end{align}
where $\psi_x(\cdot):=\psi(\tau_x+\cdot)$ and $\ph_x(\cdot):=\ph(\tau_x+\cdot)$. The last equation follows since $\tau-\tau_x$ is the first time that the state process $\ph_x$ hits $a$. That is, $\psi_x$ generates a greater payoff for the maximizer and the associated state process does not cross $x$ upwards.
Therefore, for every $x\in[a,b]$
\begin{align}\label{eq_900}
U(x)=\max\left\{0,\inf_{\pi\in\Pi}\;\sup_{\psi\in\calA_{x,\pi}} \;C(x,\pi,\psi,\tau)\right\},
\end{align}
where from now onwards $\calA_{x,\pi}$ is the restriction to absolutely continuous $\psi$'s that satisfy the conditions mentioned in (2) above.
\end{remark}

\subsection{Main results}\label{sec_2c}
We now present the main theorem, which states that the limit of the value functions of the stochastic model converge to the value function of the game. Moreover, we state an asymptotically optimal policy for the stochastic model.

For every $n\in\N$, set the stochastic control $\bpis(t) =\bpins(t)=\pi^*(\tilde W^n(t))$, $t\ge 0$, where $\pi^*$ is the function
\begin{align}
\pi^*(x)= \left\{\begin{array}{ll}\label{eq_pi}
               \displaystyle \frac{(\mu-r)e(x)}{\sigma^2(\frac{1}{2}(\frac{\mu-r}{\sigma})^2+\lambda-l(x))},
               &\ a\leq x< d,
\\ \\
                0, & \  d\le x,
              \end{array}
\right.
\end{align}
where
\begin{equation}\label{eq:defnd}
d:=b \wedge \inf\left\{y>a : \Pu- \int_a^y \frac{\lambda - l(u)+\frac{1}{2}(\frac{\mu-r}{\sigma})^2}{e(u)}du =0\right\}.
\end{equation}
and \footnote{We use the convention that $\inf\emptyset=\infty$. Also, hereafter, in case that $b=\iy$ then by the notation $(x,b]$ and $[x,b]$ mean $(x,\iy)$ and $[x,\iy)$ respectively.}
\begin{align}\notag
b:=\inf\{x\geq a : e(x)< 0\}.
\end{align}
By the definition of $b$ and since $e(\cdot)\leq M_0$ it follows that $\pi^*\in\Pi_{M_1}$ for some suitable $M_1$. 

We will show that the control $\pi^*$ is an optimal control for the minimizer in the differential game and the value function, $U$, is given by
\begin{align}
U(x)=\left\{\begin{array}{ll}\label{eq_U}
               \displaystyle \Pu- \int_a^x \frac{\lambda - l(u)+\frac{1}{2}(\frac{\mu-r}{\sigma})^2}{e(u)}du,
               &\ a\leq x< d,
\\ \\
                0, & \  d\le x.
              \end{array}
\right.
\end{align}
Since $U(x)=0$ for every $x\ge d$, the parameter $d$ is referred as the ``safe level''.
Notice that the punishment cost $\Pu$ affects $\pi^*$ and $U$ through the parameter $d$, which increases as a function of $\Pu$. Therefore, if $\Pu$ is higher, the safe level is greater. Clearly, it also affects $U$ directly linearly on $[a,d)$.

The next theorem connects between the game and the stochastic model.
\begin{theorem}[Main Result]\label{main_theorem}
Let $U^n(\cdot) := \inf_{\bpi\in\bPi}\; J^n(\cdot, \bpi)$.
For every $x \ge a$ one has, $\lim_{n\to\iy} U^n(x) = U(x)$, and moreover, $\lim_{n\to\iy}J^n(x,\bpis) = U(x)$.
\end{theorem}
The proof is given in Section \ref{sec_4}.

\section{Solution and analysis of the game}\label{sec_3}
\beginsec
In this section we provide a solution of the game.
We start by some basic properties of the value function in Section \ref{sec_3a}. In Section \ref{sec_3b} we present the HJB equation and a verification lemma. Then we derive the explicit expressions for $U$ and the optimal control. Finally, in Section \ref{sec_3e} we provide a simple control for the maximizer, which assures her the payoff $U(x)$.

\subsection{Basic properties}\label{sec_3a}

We begin by providing some basic properties that the value function satisfies. These properties are used in the verification lemma below.

\begin{lemma}\label{lem_31}
The function $U$, defined in \eqref{eq:defn-U}, satisfies the following conditions:
$\\$i. $0\leq U(x)\leq \Pu$, $x\in[a,\infty)$ and $U(a)=\Pu$.
$\\$ii. For every $x\geq b$ one has $U(x)=0$
$\\$iii. $U$ is non-increasing.
\end{lemma}
Due to part ii of the previous, in the sequel, we analyze $U$ only on the interval $[a,b)$. 

\skp
\noi{\bf Proof of Lemma \ref{lem_31}:}
i. By choosing $T=0$, the maximizer can guarantee $U\geq 0$. On the other hand, since $l(\cdot)-\lambda <0$ then clearly $U(\cdot)\leq \Pu$. By choosing $T=0$, one easily gets that $U(a)=\Pu$.

\skp
\noi
ii. We show that when  $x\ge b$, taking $\pi\equiv 0$ we can avoid ruin. Indeed, if $\pi\equiv 0$, then
$\dot\ph=-e(\ph)$, $\ph(0)=x\ge b$. In case that $x=b$ then $\ph(\cdot)\equiv b$. In case that $x>b$ then since the function $e(\cdot)$ is Lipschitz we get by Picard-Lindel\"of theorem that there is a unique $\ph\in\calC^1[0,\iy)$ that solves the ordinary differential equation that is mentioned above. One can easily verify that once $\ph$ reaches the level $b$ it remains at this level from this time onwards. Therefore,  $\ph\ge b$.

\skp
\noi
iii. 
Fix $x\in(a,\infty)$ and $y>x$ and set $\ph(0)=y$. Let $\tau_x:=\inf\{t\geq 0 : \ph(t)=x\}$. Using the dynamic programming principle along with the fact that $l(\cdot)<\lambda$, we obtain that
\begin{align}\notag
U(y)=\underset{\pi\in\Pi}{\inf}\;\;\underset{\psi\in\A_{y,\pi},T\in \R_+}{\sup}\;\left[\int_0^{T\wedge\tau_x}[-\lambda+l(\ph(t))-\frac{1}{2}\dot\psi^2(t)]dt +U(x)\right]\le U(x).
\end{align}

\noi\hfill $\Box$

\subsection{The HJB equation}\label{sec_3b}
In this section we prove that equation \eqref{eq_U} holds. We start with a verification lemma in which we provide the HJB (or rather the Isaacs) equation for the problem. Then we present a solution for this equation. Recall that by Lemma \ref{lem_31}.ii, $U(x)=0$ for $x\ge b$. Therefore, we limit ourselves to the interval $[0,b)$.

\begin{lemma}[Verification Lemma]\label{lem_verification}
Let $V:[a,b)\to[0,\Pu]$, with $V(a)=\Pu$, be a non-increasing, continuous function that is differentiable on $(a,\beta)$, where $\beta:=b\wedge\inf\{x>a: V(x)=0\}$. Assume that the following conditions hold:

\noi
(i) For every $x\in[a,\beta)$ one has
\begin{align}[\text{The HJB equation}]\qquad\label{eq_901}
\inf_{p\in\R}\;\sup_{\theta\in\R}\;\left\{V'(x)(-e(x)+(\mu-r)p+\sigma p\theta)-\lambda+l(x)-\frac{1}{2}\theta^2\right\}=0;
\end{align}

\noi
(ii) Let $P(x)=-\frac{\mu-r}{\sigma^2 V'(x)}$. Then for every $\psi \in \mathcal{AC}_0$ and every $x \in [a,\beta)$, there exists a unique solution to \eqref{eq_250} when we replace $\pi$ with $P$.


\skp
\noi
(iii) Let $\Theta(p, x)=\sigma p V'(x)$. Then for every $\pi\in\Pi$, and $a\le x<\beta$, there exists a unique solution to 
\begin{align}\label{new100}
\dot\ph(t)&=-e(\ph(t))+(\mu-r)\pi(\ph(t)) + \sigma\pi(\ph(t))\Theta(\pi(\ph(t)),\ph(t)),\;\;t\in[0,\tau_a],\\\notag
\ph(0)&=x
\end{align}
such that $\int_0^t\Theta(\pi(\ph(u)),\ph(u))du \in\calA_{x,\pi}$ (see Remark \ref{rem_21}), where $\tau_a:=\inf\{t\ge 0:\ph(t)=a\}$.

\skp
\noi
Then $U=V$ on $[a,b)$. Moreover,  the function $P$ is an optimal feedback control.

\end{lemma}
Notice that we defined the HJB equation only on the interval $[a,\beta)$. This structure follows since for every $x$ for which $U(x)=0$, under optimality of both players, the time part of the maximizer's control equals zero and the game is terminated immediately. 

\skp\noi
\textbf{Proof of Lemma \ref{lem_verification}:}
1.  As a first step we will prove that for every $x\in[a,\beta)$ one has $V(x)\ge U(x)$. As a result, if $\beta< b$ then $0=V(\beta)\ge U(\beta)\ge 0$, where the last inequality follows by the first assertion of Lemma \ref{lem_31}. Since $V\ge 0$ and $U$ is non-increasing we get that $V\ge U$ on $[a,b)$.

Fix $x\in[a,\beta)$. Set the control $\pi^*=P$. Also, fix a control $\psi\in\A_{x,\pi^*}$ and denote by $\ph^*$ the state process associated with $\pi^*$ and $\psi$. Recall that by the definition of $\A_{x,\pi^*}$, for every $t> 0$, one has $\ph^*(t)<\ph^*(0)=x$ and that $\tau^*_a<\iy$, where $\tau^*_a$ is the first time that $\ph^*$ reaches $a$. Recalling moreover that $x<\beta$ we get that for every $t\ge 0$, $\ph^*(t)<\beta$. Since $V$ is differentiable on $[a,\beta)$ we can apply the chain rule to $V$ and get
\begin{align}\notag
V(\ph^*(\tau^*_a))-V(\ph^*(0)) = \int_0^{\tau^*_a}V'(\ph^*(t))[-e(\ph^*(t))+(\mu-r)\pi^*(\ph^*(t)) + \sigma\pi^*(\ph^*(t))\dot\psi(t)]dt.
\end{align}
Using again the inequality  $\ph^*(t)<\beta$ we get by conditions (i) and (ii) that
\begin{align}\notag
0&=\sup_{\theta\in\R}\;\left\{V'(\ph^*(t))[-e(\ph^*(t))+(\mu-r)\pi^*(\ph^*(t))+\sigma\pi^*(\ph^*(t))\theta]-\lambda+l(\ph^*(t))-\frac{1}{2}\theta^2\right\}\\\notag
&\ge V'(\ph^*(t))[-e(\ph^*(t))+(\mu-r)\pi^*(\ph^*(t))+\sigma\pi^*(\ph^*(t))\dot\psi(t)]-\lambda+l(\ph^*(t))-\frac{1}{2}\dot\psi^2(t).
\end{align}
So we have
\begin{equation}\label{eq:gtopt}
V(x) \ge\int_0^{\tau^*_a}[-\lambda+l(\ph^*(t))-\frac{1}{2}\dot\psi^2(t)]dt+\Pu.
\end{equation}
By taking first $\sup_{\psi\in\A_{x,\pi}}$ and then $\inf_{\pi\in\Pi}$ on both sides, we get that
$$V(x)\ge \inf_{\pi\in\Pi}\;\sup_{\psi\in\A_{x,\pi}}\; C(x,\pi,\psi,\tau^*_a).$$
By the above and recalling that $V\ge 0$ we get by \eqref{eq_900} that $V(x)\ge U(x)$.

\skp\noi
2. By using the assumption that $V$ is non-increasing and that $U\ge 0$ it follows that it is sufficient to prove that $V\le U$ on $x\in[a,\beta)$. Fix $x\in[a,\beta)$. For every $\pi\in\Pi$ denote $\Theta(\pi(\ph^*(t)),\ph^*(t))$ by $\dot\psi^*(t)$, where $\ph^*$ solves \eqref{new100}.
By (iii) $\ph^*$ reaches $a$ in a finite time, that is, $\tau_a<\iy$. 
Using conditions (i) and (ii) we get that
\begin{align}\label{eq_906}
0&=\inf_{p\in\R}\;\left\{V'(\ph^*(t))[-e(\ph^*(t))+(\mu-r)p+\sigma p\Theta(p,\ph^*(t))]-\lambda+l(\ph^*(t))-\frac{1}{2}\Theta^2(p,\ph^*(t))\right\}\\\notag
&\le V'(\ph^*(t))[-e(\ph^*(t))+(\mu-r)\pi(\ph^*(t))+\sigma\pi(\ph^*(t))\dot\psi^*(t)]-\lambda+l(\ph^*(t))-\frac{1}{2}(\dot\psi^*)^2(t).
\end{align}
Recalling that 
\begin{align}\notag
\dot\ph^*(t) = -e(\ph^*(t))+(\mu-r)\pi(\ph^*(t))+\sigma\pi(\ph^*(t))\dot\psi^*(t), 
\end{align}
we get that
\begin{align}\notag
V(x) \le\int_0^{\tau_a}[-\lambda+l(\ph^*(t))-\frac{1}{2}(\dot\psi^*)^2(t)]dt+\Pu.
\end{align}
By taking $\sup_{\psi\in\A_{x,\pi}}$ first and then $\inf_{\pi\in\Pi}$ on both sides, we get that
$$V(x)\le \inf_{\pi\in\Pi}\;\sup_{\psi\in\A_{x,\pi}}\; C(x,\pi,\psi,\tau) .$$
By the definition of $\beta$ and since $x<\beta$ it follows that $V(x)>0$. Therefore, $$\inf_{\pi\in\Pi}\;\sup_{\psi\in\A_{x,\pi}}\; C(x,\pi,\psi,\tau)>0$$ and by \eqref{eq_900} it follows that $U(x) = \inf_{\pi\in\Pi}\;\sup_{\psi\in\A_{x,\pi}}\; C(x,\pi,\psi,\tau)$. Therefore, $V(x)\le U(x)$. Now, the optimality of the feedback control $P$ follows from \eqref{eq:gtopt}. \hfill $\Box$ \\

We now use the verification lemma in order to provide an explicit expression for the value function $U$.
\begin{proposition}\label{prop_U}
Let
\begin{align}
V(x)=\left\{\begin{array}{ll}\label{eq_909b}
               \displaystyle \Pu- \int_a^x \frac{\lambda - l(u)+\frac{1}{2}(\frac{\mu-r}{\sigma})^2}{e(u)}du,
               &\ a\leq x< d,
\\ \\
                0, & \  d\le x,
              \end{array}
\right.
\end{align}
where $d$ is defined in \eqref{eq:defnd}. Then $V=U$, defined in \eqref{eq:defn-U}. Moreover, the control $\pi^*$ given in \eqref{eq_pi} is optimal.
\end{proposition}
\noi
\textbf{Proof:} Recall that by the second assertion of Lemma \ref{lem_31}, for every $x\ge b$, $U(x)=0$. We now prove that \eqref{eq_909b} holds for $x\in[a,b)$ using Lemma \ref{lem_verification}. Notice that the parameter $\beta$ that appears in the verification lemma is actually $d$ for our particular function $V$. First, it is easy to check that $V$ satisfies \eqref{eq_901}. Next, since $P:[a,d] \to \mathbb{R}$ given by
\[
P(x)=-\frac{\mu-r}{\sigma^2 V'(x)}=\frac{(\mu-r)e(x)}{\sigma^2(\frac{1}{2}(\frac{\mu-r}{\sigma})^2+\lambda-l(x))}
\]
is (locally) Lipschitz, requirement ii. of the verification lemma also holds.

Now we will verify the third condition in the verification lemma. 
Notice, that \eqref{new100} admits a unique solution on the time interval $[0,\tau_a\wedge\tau_d]$, where $\tau_d$ is the first time that $\ph$ hits $d$. The equation in \eqref{new100} 
admits a unique solution on $[0,\tau_a\wedge\tau_d]$ since (a) on this time interval $\ph\in[a,d]$ and the functions $e,\pi$, and $V'$ are bounded on the interval $[a,d]$ and (b) $e$ is Lipschitz, and $\pi^2$ and $V'$ are locally Lipschitz. 

Next, we will show that the above argument can be upgraded to the interval $[0,\tau_a]$.
Observe that the inequality in \eqref{eq_906} holds with $\dot\psi(\cdot):=\Theta(\pi(\ph(\cdot)),\ph(\cdot))$ and $\ph(t)$ replacing $\dot\psi^*(t)$ and $\ph^*(t)$. Moreover, since $V$ is non-increasing, we get that $V'(\ph(t))\le 0$ and together with \eqref{eq_906} in this case, actually $V'(\ph(t))<0$ on $t\in[0,\tau_a\wedge\tau_d]$. Hence, 
\begin{align}\notag
\dot\ph(t) = -e(\ph(t))+(\mu-r)\pi(\ph(t))+\sigma\pi(\ph(t))\dot\psi(t) \le \frac{\lambda-l(\ph(t))+\frac{1}{2}(\dot\psi)^2(t)}{V'(\ph(t))}<0.
\end{align}
Thus, $\ph$ does not cross $\ph(0)$ upwards and therefore, $\tau_d=\iy$, which implies \eqref{new100}. 

As a final step we will show that  $\ph$ hits $a$ in a finite time as a result of which we will obtain that
we get that $\psi\in\A_{x,\pi}$. Assume to the contrary that $\tau_a=\iy$, then by using again \eqref{eq_906} in our case, we conclude that for every $s>0$ one has
\begin{align}\notag
V(\ph(0))-V(\ph(s)) \le\int_0^{s}[-\lambda+l(\ph(t))-\frac{1}{2}(\dot\psi)^2(t)]dt.
\end{align}
Since $l(\cdot)<\lambda$ we get that the r.h.s.~of the above goes to $-\iy$ when $s\to\iy$, which contradicts the fact that $V$ is bounded.

The optimality of $\pi^*$ follows by the verification lemma and since $U=V$.

\hfill $\Box$

\subsection{Saddle point property}\label{sec_3e}

Here, we provide a control (for the maximizer) that is independent of $\pi$ and that assures her the payoff $U(x)$. The simplicity of this control will be crucial in Section \ref{sec_4a}. Set
\begin{align}\label{eq_365}
\psi^*(t) = -\frac{\mu-r}{\sigma}t,\quad t\ge 0,
\end{align}
and let $T^*=\tau$ in case $x<d$ and $T^*=0$ otherwise. Notice that $\psi^*$ is independent of the control $\pi$. Moreover, notice that under $\psi^*$ the state process satisfies $\dot\ph = -e(\ph)$ and therefore $\ph$ and $T^*$ are also independent of the choice of $\pi$.

\begin{proposition}\label{prop_34}
For every $x\in[a,\iy)$ one has $U(x)=\underset{\pi\in\Pi}{\inf} \;C(x,\pi,\psi^*,T^*) $.
Moreover,
\begin{align}\label{eq_366}
T^*\leq \frac{\Pu-U(x)}{\lambda-l(a)+\frac{1}{2}(\frac{\mu-r}{\sigma})^2}.
\end{align}
\end{proposition}
\noi
{\bf Proof:}
For $x\geq d$ one has $U(x)=0$ and by definition $T^*=0$, so that $\underset{\pi}{\inf} \;C(x,\pi,\psi^*,T^*)=0$.
Set $x<d$. 

Under $\psi^*$, the state process is independent of $\pi$. Hence, for every $\pi\in\Pi$
\begin{align}\notag
C(x,\pi,\psi^*,T^*)&=\int_0^{T^*}\left[-\lambda+l(\ph(t))-\frac{1}{2}(\dot\psi^*)^2(t)\right]dt + \Pu\\\notag
                             &=-\int_a^x\frac{\lambda-l(u)+\frac{1}{2}(\frac{\mu-r}{\sigma})^2}{e(u)}du + \Pu\\\notag
                             & =U(x),
\end{align}
where the second equality follows by the change of variables $u=\ph(t)$, and the last equality follows by Proposition \ref{prop_U}. Hence, the first part of the theorem is proved.

Since $l$ is non-increasing we get by \eqref{eq_900} that
\begin{align}\notag
0<U(x) &=    \int_0^{T^*}\left[-\lambda+l(\ph(t))-\frac{1}{2}\left(\frac{\mu-r}{\sigma}\right)^2\right]dt + \Pu\\\notag
     &\leq \int_0^{T^*}\left[-\lambda+l(a)-\frac{1}{2}\left(\frac{\mu-r}{\sigma}\right)^2\right]dt + \Pu\\\notag
     &=    -T^*\left[\lambda-l(a)+\frac{1}{2}\left(\frac{\mu-r}{\sigma}\right)^2\right] + \Pu
\end{align}
and \eqref{eq_366} follows.
\hfill $\Box$

\section{Proof of Theorem \ref{main_theorem}}\label{sec_4}
\beginsec
The proof of Theorem \ref{main_theorem} follows by some measure changing arguments and also influenced by Varadhan's lemma. We show in two separate theorems that $U(x)$ is a lower (resp., upper) bound to $\lim\inf_{n\to \infty} U^n(x)$ (resp., $\lim\sup_{n\to\infty} U^n(x)$), where $U^n(\cdot) := \inf_{\bpi\in\bPi}\; J^n(\cdot, \bpi)$. Moreover, we show that the policy $\pi^*$ is asymptotically optimal.



\subsection{Lower bound}\label{sec_4a}
\begin{theorem}\label{thm_45}
For every $x\ge a$ one has $\lim\inf_{n\to\infty} U^n(x)\ge U(x)$.
\end{theorem}
\noi
\textbf{Proof:}
Recall that $U(x)=0$ for every $x\geq d$. Since $l\ge 0$ and also $\rho>0$ it follows from Jensen's inequality that for any sequence of policies $\{\bpin\}_n$, we have  
\begin{align}\notag
\frac{1}{n}\ln \E\left[e^{n\left(\int_0^{\tau^n_d\wedge\tau^n_a}l(\tilde W^n(t))dt+\Pu 1_{\{\tau^n_a\leq \tau^n_d\}}\right)}\right]\ge 0.
\end{align}

Fix $x\in(a,d)$ and fix an arbitrary sequence of policies $\{\bpin\}_n\subseteq\Pi$. We show that for every $\eps>0$ there is $N>0$ such that for every $n>N$ one has $J^n(x,\bpin)\geq U(x)-w_0(\eps)$, where $w_0(\eps)\to0$ as $\eps\to0$.

We start with some preliminaries. Let $\psi^*$ be the function from \eqref{eq_365}
and let
\begin{align}\label{eq_409}
\dot\ph^*(t)=\left\{\begin{array}{ll}
               \displaystyle  -e(\ph^*(t)),
               &\ 0\leq t\leq T^*,
\\ \\
                -e(a), & T^*\leq t  ,
              \end{array}
\right.
\end{align}
with $\ph^*(0)=x$, where $T^*$ is the first time that $\ph^*$ hits $a$, which is finite thanks to \eqref{prop_34}.
Note that up to time $T^*$, $\ph^*$ is the state process of the differential game associated with $\psi^*$ and any control $\pi\in\Pi$.

Let us fix $\eps_1>0$. Since $l$ is Lipschitz then there exists $\gamma_1>0$ such that for every $y,z\in[a,\iy)$
\begin{align}\label{eq_410}
|y-z|<\gamma_1\qquad \text{implies}\qquad |l(y)-l(z)|<\eps_1.
\end{align}
Moreover, since $\ph^*$ is continuous and for $t>T^*$, $\dot\ph^*(t)<0$ it follows that one may choose $\gamma_1$ such that
\begin{align}\label{eq_411}
|\ph-\ph^*|_{T^*+2\eps_1} \leq \gamma_1\qquad \text{implies}\qquad |\tau_a[\ph]-T^*|\leq \eps_1,
\end{align}
where $\tau_a[\ph]:=\inf\{t\geq 0: \ph(t)=a\}$. Indeed, recall that $\ph(0)=x\in(a,d)$. Now, since $e(\cdot)$ is positive on $[a,d)$ we get from \eqref{eq_409} that the state process $\ph^*$ is strictly decreasing on $[0,T^*+2\eps_1]$, touching $a$ only at $T^*$ and continuing to decrease on $[T^*,T^*+2\eps_1]$.

Define the probability measure $Q^*=Q^{*,n}$ on $(\Omega,\calF_{T^*+2\eps_1})$ by
\begin{align}\notag
\frac{dQ^*}{d\PP}(t) = e^{-\sqrt{n}\int_0^t\dot\psi^*(s)dB(s) - \frac{n}{2}\int_0^t(\dot\psi^*)^2(s)ds },\quad t\in[0,T^*+2\eps_1].
 \end{align} 
Then under $Q^*$, $ B^*(t)=B^{*,n}(t):=B(t)+\sqrt{n}\frac{\mu-r}{\sigma}t$, $t\in[0,T^*+2\eps_1]$ is a standard Brownian motion and 
\begin{align}\notag
d\tilde W^n(t) &= -e(\tilde W^n(t))+\frac{1}{\sqrt{n}}\sigma\bpin(\tilde W^n(t))dB^*(t),\quad t\in[0,T^*+2\eps_1].\\\notag
\end{align}
Now, since that $|\bpin(t)|\le M_1$ then by Gronwall's inequality and Doob's martingale inequality we get that there is a constant $C_1>0$ that depends on the Lipschitz constant of $e(\cdot)$ such that 
\begin{align}\label{new1}
Q^*\left(\left(\calE^n\right)^c\right)\le\frac{C_1 M_1^2}{n\gamma_1^2},
\end{align}
where 
$$\calE^n:=\left\{\omega:\left|\tilde W^n(\cdot,\omega)-\ph^*(\cdot,\omega)\right|_{T^*+2\eps_1}\le \gamma_1\right\}.$$
Set $N=N(\eps_1,\gamma_1,M_1,C_1)$ such that 
\begin{align}\label{new2}
N>\max\left\{\frac{-\ln(\eps_1)}{\eps_1},\frac{C_1M_1^2}{\eps_1\gamma_1^2},\frac{C_1M_1^2(T^*+2\eps_1)\left(\lambda+\tfrac{1}{2}\left(\tfrac{\mu-r}{\sigma}\right)^2\right)}{\eps_1\gamma_1^2}\right\}.
\end{align}

We are now ready to bound from below $J^n(x,\bpin)$.
Fix $n>N$ then
\begin{align}\notag
&\frac{1}{n}\ln\;\E\left[\int_0^{\infty}e^{-\lambda n t}e^{n\left(\int_0^{\tau^n_a\wedge t}l(\tilde W^n(s))ds + \Pu 1_{\{\tau^n_a\leq t\}}\right)}dt\right]\\\notag
&\quad\geq
\frac{1}{n}\ln\;\E\left[\int_{T^*+\eps_1}^{T^*+2\eps_1}e^{-\lambda n t}e^{n\left(\int_0^{\tau^n_a\wedge t}l(\tilde W^n(s))ds + \Pu 1_{\{\tau^n_a\leq t\}}\right)}dt\right]\\\notag
&\quad\geq
\frac{1}{n}\ln\;\E\left[\eps_1 e^{-\lambda n (T^*+2\eps_1)}e^{n\left(\int_0^{\tau^n_a\wedge (T^*+\eps_1)}l(\tilde W^n(s))ds + \Pu 1_{\{\tau^n_a\leq T^*+\eps_1\}}\right)}\right]\\\notag
&\quad=\frac{1}{n}\ln\;\E^{Q^*}\left[\eps_1 e^{-\lambda n (T^*+2\eps_1)}e^{n\left(\int_0^{\tau^n_a\wedge (T^*+\eps_1)}l(\tilde W^n(s))ds + \Pu 1_{\{\tau^n_a\leq T^*+\eps_1\}}\right)}\frac{d\PP}{dQ^*}(T^*+2\eps_1)\right]\\\notag
&\quad\ge\;\E^{Q^*}\left[-\lambda (T^*+2\eps_1)+\int_0^{\tau^n_a\wedge (T^*+\eps_1)}l(\tilde W^n(s))ds + \Pu 1_{\{\tau^n_a\leq T^*+\eps_1\}}-\frac{1}{2}\int_0^{T^*+2\eps_1}(\dot\psi^*)^2(s)ds\right]\\\notag&\qquad-\eps_1\\\notag
&\quad\ge\;\E^{Q^*}\left[\left(-\lambda (T^*+2\eps_1)+\int_0^{\tau^n_a\wedge (T^*+\eps_1)}l(\tilde W^n(s))ds + \Pu 1_{\{\tau^n_a\leq T^*+\eps_1\}} -\frac{1}{2}\int_0^{T^*+2\eps_1}(\dot\psi^*)^2(s)ds\right)1_{\calE^n}\right]\\\notag&\qquad-2\eps_1\\\notag
&\quad\ge\;\E^{Q^*}\left[\left(-\lambda (T^*+2\eps_1)+\int_0^{T^*-\eps_1}[l(\ph^*(s))-\eps_1]ds + \Pu -\frac{1}{2}\int_0^{T^*+2\eps_1}(\dot\psi^*)^2(s)ds\right)1_{\calE^n}\right]-2\eps_1\\\notag
&\quad=\left(-\lambda (T^*+2\eps_1)+\int_0^{T^*-\eps_1}[l(\ph^*(s))-\eps_1]ds + \Pu -\frac{1}{2}\int_0^{T^*+2\eps_1}(\dot\psi^*)^2(s)ds\right)Q^*(\calE^n)-2\eps_1\\\notag
&\quad=-\lambda T^*+\int_0^{T^*}l(\ph^*(s))ds-\frac{1}{2}\int_0^{T^*}(\dot\psi^*)^2(s)ds+\Pu+w_0(\eps_1)\\\notag
&\quad=U(x)+w_0(\eps_1),\\\notag
\end{align}
where 
\begin{align}\notag
w_0(\eps_1)&=\left(-2\lambda\eps_1-\int_{T^*-\eps_1}^{T^*}l(\ph^*(s))ds-\eps_1(T^*-\eps_1)-\frac{1}{2}\int_{T^*}^{T^*+2\eps_1}(\dot\psi^*)^2(s)ds\right)Q^*(\calE^n)-2\eps_1\\\notag
&\quad-\left(
-\lambda T^*+\int_0^{T^*}l(\ph^*(s))ds-\frac{1}{2}\int_0^{T^*}(\dot\psi^*)^2(s)ds+\Pu
\right)(1-Q^*(\calE^n)).
\end{align}
The first three relations are easy to check. The forth relation follows by Jensen's inequality and by \eqref{new1}. The fifth relation follows since $l(\cdot)\ge 0$ and by \eqref{new1} and \eqref{new2}. The sixth relation follows by \eqref{eq_410} and \eqref{eq_411}. The seventh relation follows since all the terms inside the expectation besides the indicator are deterministic. Finally, the last relation follows by Proposition \ref{prop_34}. By \eqref{new1} and \eqref{new2} and recalling that $-\lambda T^*+\int_0^{T^*}l(\ph^*(s))ds-\frac{1}{2}\int_0^{T^*}(\dot\psi^*)^2(s)ds+\Pu=U(x)\ge 0$ we get that $w_0(\eps_1)\to 0$ as $\eps_1\to 0$.

\hfill$\Box$

\subsection{Asymptotically optimal policy}
In this section we show that the optimal policy in the game, which was defined in \eqref{eq_pi}
is an asymptotically optimal policy in the stochastic model. We start with a technical lemma that serves us in the proof of the preceding theorem. The lemma provides an upper bound for the discounted cost by an alternative cost that is defined through a new measure $Q^n$. Set $T:=\Pu /(\lambda-l(a))$. Both $Q^n$ and $T$ will play important roles during the proof of the theorem.

\begin{lemma}\label{lem_45}
For every $n\in\N$, there exists a probability measure $Q^n\sim\PP$ on $[0,T]$, for which
\begin{align}\label{eq_426}
\lim&\sup\frac{1}{n}\ln \E\left[\int_0^\infty e^{-\lambda nt}e^{n\left(\int_0^{\tau^n_a\wedge t}l(\tilde W^n(s))ds +\Pu 1_{\{\tau^n_a\leq t\}}\right)}dt\right]\\\notag
&\leq \lim\sup \E^{Q^n}\left[
\underset{0\leq t\leq T}{\sup}\left(\int_0^{\tau^n_a\wedge t}[l(\tilde W^n(s))-\lambda ]ds +\Pu 1_{\{\tau^n_a\leq t\}} \right)\right]-\frac{1}{n}\calH(Q^n\|\PP),
\end{align}
where $\calH(Q^n\|\PP):=E^{Q^n}\left[\ln\Big(\tfrac{dQ^n}{d\PP}\Big)\right]$ is the relative entropy of $Q^n$ w.r.t.~$\PP$.
Also, there is $N_1\in\N$ such that for every $n>N_1$ one has
\begin{align}\label{new7}
 \frac{1}{n}\calH(Q^n\|\PP)\le 2\Pu.
\end{align} 
Moreover, for every $n\in\N$, there exists an adapted process $(\psi(t))_{0\le t\le T}$ such that $Q^n$-almost surely, $\psi^n(\cdot,\omega)\in\calA\calC_0[0,T]$ and $\int_0^{T}(\dot\psi^n)^2(s,\omega)<\iy$, and
\begin{align}\label{cor1}
\frac{1}{n}\calH(Q^n\|\PP) =E^{Q^n}\left[ \frac{1}{2}\int_0^{T}(\dot\psi^n)^2(s)ds\right].
\end{align}
\end{lemma}
\noi
\textbf{Proof:}
First, notice that 
\begin{align}\notag
&\E\left[\int_0^\infty e^{-\lambda nt}e^{n\left(\int_0^{\tau^n_a\wedge t}l(\tilde W^n(s))ds +\Pu 1_{\{\tau^n_a\leq t\}}\right)}dt\right]\\\notag
 &\quad\leq
\E\left[\int_0^T e^{-\lambda nt}e^{n\left(\int_0^{\tau^n_a\wedge t}l(\tilde W^n(s))ds +\Pu 1_{\{\tau^n_a\leq t\}}\right)}dt+e^{n\Pu} \int_T^\infty e^{nt\left(l(a)-\lambda \right)}dt\right]\\\notag
 &\quad=
\E\left[\int_0^T e^{-\lambda nt}e^{n\left(\int_0^{\tau^n_a\wedge t}l(\tilde W^n(s))ds +\Pu 1_{\{\tau^n_a\leq t\}}\right)}dt+  \frac{1}{n(\lambda-l(a))}e^{n\left(\Pu+T\left(l(a)-\lambda\right) \right)}\right]\\\notag
 &\quad\leq
\E\left[\int_0^T e^{n\left(\int_0^{\tau^n_a\wedge t}[l(\tilde W^n(s))-\lambda ]ds +\Pu 1_{\{\tau^n_a\leq t\}}\right)}dt+  \frac{1}{n(\lambda-l(a))}\right]\\\notag
 &\quad\leq
T\E\left[ e^{n\underset{0\leq t\leq T}{\sup}\left(\int_0^{\tau^n_a\wedge t}[l(\tilde W^n(s))-\lambda ]ds +\Pu 1_{\{\tau^n_a\leq t\}}\right)}+  \frac{1}{\rho n}\right].
\end{align}
The first inequality follows since $l$ is non-increasing and since by eliminating the indicator in the second exponent we only increase the cost. The second inequality follows by the choice of $T$ and since $-\lambda n t\le -\lambda n(\tau^n_a\wedge t)$. The other relations are easy to see.

Now, for every $n\in\N$ let $Q^n$ be the measure that satisfies
\begin{align}\label{new4}
&\frac{1}{n}\ln \left( \E\left[e^{n\underset{0\leq t\leq T}{\sup}\left(\int_0^{\tau^n_a\wedge t}[l(\tilde W^n(s))-\lambda ]ds +\Pu 1_{\{\tau^n_a\leq t\}}\right)}\right]\right)\\\notag
&\quad=\E^{Q^n}\left[
\underset{0\leq t\leq T}{\sup}\left(\int_0^{\tau^n_a\wedge t}[l(\tilde W^n(s))-\lambda ]ds +\Pu 1_{\{\tau^n_a\leq t\}} \right)
\right]-\frac{1}{n}\calH(Q^n\|\PP).
\end{align}
The existence of the measure $Q^n$ with the above representation is justified by the following argument: It follows from Jensen's inequality that for any measure $Q\sim \PP$ and 
\[
f=
n \underset{0\leq t\leq T}{\sup}\left(\int_0^{\tau^n_a\wedge t}[l(\tilde W^n(s))-\lambda ]ds +\Pu 1_{\{\tau^n_a\leq t\}} \right)
\]
 one has 
\begin{align}\notag
\ln\left(\E\left[e^f\right]\right) = \ln\left(\E^Q\left[e^f\cdot\frac{dP}{dQ}\right]\right) \ge \E^Q\left[f\right]-\calH(Q\|\PP),  
\end{align}
where equality holds for the measure $Q$ that satisfies
\begin{align}\label{new11}
\frac{dQ}{dP}=\frac{e^f}{\E\left[e^f\right]}.  
\end{align}

By point \emph{i} of Lemma \ref{lem_31}, Theorem \ref{thm_45}, and \eqref{new4} we get that
\begin{align}\notag
0&\le U(x)\le \liminf_{n\to\infty} J^n(x,\bpis)  \leq \liminf_{n \to \infty}\left(\frac{1}{n}\E^{Q^n}\left[f\right]-\frac{1}{n}\calH(Q^n\|\PP)\right)\le\rho- \limsup_{n \to \infty}\left(\frac{1}{n}\calH(Q^n\|\PP)\right),
\end{align}
where the last inequality follows since $\lambda>l(\cdot)$. Hence, \eqref{new7} holds.

We now turn to the last part of the lemma. Since the r.h.s.~of \eqref{new11} conditioned on $\calF_t$ is a positive $P$-martingale, then it can be expressed as an exponential martingale. That is, there is a predictable and square integrable process $(u^n(t))_{0\le t\le T}$ such that 
\begin{align}\notag
\frac{dQ^n}{d\PP}(t) = e^{\sqrt{n}\int_0^t u^n(s)dB(s) - \frac{n}{2}\int_0^t(u^n)^2(s)ds },\quad t\in[0,T].
\end{align}
 Now, for $Q^n$-almost every $\omega$ define $\psi^n(\cdot,\omega)$ as the Lesbegue integral of $u^n(\cdot,\omega)$.

Finally, notice that under $Q^n$, $B(t)- \sqrt{n}u^n(t)$, $t\in[0,T]$ is a Brownian motion and therefore \eqref{cor1} holds.
 
\hfill$\Box$

\begin{theorem}\label{thm_46}
For every $x\ge a$ one has $\lim\sup_{n\to\infty} J^n(x,\bpis)\le U(x)$, where $\bpis$ is defined on \eqref{eq_pi}.
\end{theorem}
\noi
\textbf{Proof:}
Fix $x\ge a$. Notice that by \eqref{eq_426} it is sufficient to bound the $\limsup$ of its r.h.s.~by $U(x)$.

During the proof we will make use of the state process and the wealth process under the control $\bpis$, which was defined via the function $\pi^*$, see \eqref{eq_pi}. Consider $Q^n$ and $\psi^n$ from Lemma \ref{lem_45}. Under $Q^n$
\begin{align}\label{refe4a}
d\tilde W^n(t) &= \left(-e(\tilde W^n(t))+(\mu-r)\pi^*(\tilde W^n(t))+\sigma\pi^*(\tilde W^n(t))\dot\psi^n(t)\right)dt+\frac{1}{\sqrt{n}}\sigma\pi^*(\tilde W^n(t))dB^n(t),
\end{align}
$t\in(0,T]$ and $\tilde W^n(0)=x$, where $B^n(t)=B(t)-\sqrt{n}\psi^n(t)$ is a standard Brownian motion under $Q^n$. Also, set
\begin{align}\label{refe4b}
\dot\ph^n(t)&=-e(\ph^n(t))+(\mu-r)\pi^*(\ph^n(t)) + \sigma\pi^*(\ph^n(t))\dot\psi^n(t),\;\;t\in[0,T],\\\notag
\ph^n(0)&=x.
\end{align}
For any $\delta>0$, we define
\begin{align}\label{refe6}
A_n(\delta)=\{\omega : |\tilde W^n-\ph^n|_{\tau_a[\tilde W^n]\wedge T}\le \delta\},
\end{align}
where $\tau_a[h]:=\inf\{t\in[0,T]: h(t)\le a\}$ with the convention that $\inf \emptyset=\infty$.
Let us write
\begin{align}\label{refe1}
&\E^{Q^n}\left[\sup_{0\le t\le T} \int_0^{\tau_a[\tilde W^n]\wedge t}[l(\tilde W^n(s))-\lambda-\frac{1}{2}(\dot\psi^n(s))^2]ds+\rho 1_{\{{\tau_a[\tilde W^n]\le t\}}}\right]\\\notag
&\quad=\E^{Q^n}\left[\left(\sup_{0\le t\le T} \int_0^{\tau_a[\tilde W^n]\wedge t}[l(\tilde W^n(s))-\lambda-\frac{1}{2}(\dot\psi^n(s))^2]ds+\rho 1_{\{{\tau_a[\tilde W^n]\le t\}}}\right)1_{A_n(\delta)}\right]\\\notag
&\qquad+\E^{Q^n}\left[\left(\sup_{0\le t\le T} \int_0^{\tau_a[\tilde W^n]\wedge t}[l(\tilde W^n(s))-\lambda-\frac{1}{2}(\dot\psi^n(s))^2]ds+\rho 1_{\{{\tau_a[\tilde W^n]\le t\}}}\right)1_{(A_n(\delta))^c}\right].
\end{align}
For $\omega\in A_n(\delta)$, we have $\tau_a[\tilde W^n]\ge\tau_{a+\delta}[\ph^n]$. Then there is a constant  $c_1>0$ that depends on the Lipschitz constant of $l(\cdot)$ such that on $A_n(\delta)$ and for every $t\in[0,T]$,
\begin{align}\notag
 &\int_0^{\tau_a[\tilde W^n]\wedge t}[l(\tilde W^n(s))-\lambda-\frac{1}{2}(\dot\psi^n(s))^2]ds+\rho 1_{\{{\tau_a[\tilde W^n]\le t\}}}\\\notag
 &\quad\le \int_0^{\tau_a[\tilde W^n]\wedge t}[l(\ph^n(s))-\lambda-\frac{1}{2}(\dot\psi^n(s))^2]ds+\rho 1_{\{{\tau_a[\tilde W^n]\le t\}}}+c_1\delta\\\notag
 &\quad\le \int_0^{\tau_{a+\delta}[\ph^n]\wedge t}[l(\ph^n(s))-\lambda-\frac{1}{2}(\dot\psi^n(s))^2]ds+\rho 1_{\{{\tau_{a+\delta}[\ph^n]\le t\}}}+c_1\delta\\\notag
 &\quad\le U_{a+\delta}(x)+c_1\delta.
\end{align}
Here we denote by $U_{a+\delta}(x)$ the function defined by \eqref{eq:defn-U} with $a$ replaced by $a+\delta$. 
The last inequality follows since the optimal control $\pi^*$ defined by \eqref{eq_pi} is, according to
its explicit form, the optimal control also for the differential game with $a$
replaced by $a+\delta$.
Therefore,
\begin{align}\notag
&\E^{Q^n}\left[\left(\sup_{0\le t\le T} \int_0^{\tau_a[\tilde W^n]\wedge t}[l(\tilde W^n(s))-\lambda-\frac{1}{2}(\dot\psi^n(s))^2]ds+\rho 1_{\{{\tau_a[\tilde W^n]\le t\}}}\right)1_{A_n(\delta)}\right]\\\notag
&\quad\le(U_{a+\delta}(x)+c_1\delta)Q^n(A_n(\delta)).
\end{align}
On the other hand,
\begin{align}\notag
&\E^{Q^n}\left[\left(\sup_{0\le t\le T} \int_0^{\tau_a[\tilde W^n]\wedge t}[l(\tilde W^n(s))-\lambda-\frac{1}{2}(\dot\psi^n(s))^2]ds+\rho 1_{\{{\tau_a[\tilde W^n]\le t\}}}\right)1_{(A_n(\delta))^c}\right]\\\notag
&\quad\le\rho Q^n((A_n(\delta))^c).
\end{align}
Plugging the last two inequalities into
 \eqref{refe1} we conclude
\begin{align}\notag
&\E^{Q^n}\left[\sup_{0\le t\le T} \int_0^{\tau_a[\tilde W^n]\wedge t}[l(\tilde W^n(s))-\lambda-\frac{1}{2}(\dot\psi^n(s))^2]ds+\rho 1_{\{{\tau_a[\tilde W^n]\le t\}}}\right]\\\notag
&\quad\le  U_{a+\delta}(x)+c_1\delta+ (\rho-U_{a+\delta}(x)-c_1\delta) Q^n((A_n(\delta))^c).
\end{align}
In the following, we shall show that 
\begin{align}\label{refe3}
\lim_{n\to\iy} Q^n((A_n(\delta))^c)=0,
\end{align}
from which it follows that 
\begin{equation}\label{eq:asfeqn}
\begin{split}
&\limsup\;\E^{Q^n}\left[\sup_{0\le t\le T} \int_0^{\tau_a[\tilde W^n]\wedge t}[l(\tilde W^n(s))-\lambda-\frac{1}{2}(\dot\psi^n(s))^2]ds+\rho 1_{\{{\tau_a[\tilde W^n]\le t\}}}\right]\\
&\quad\le U_{a+\delta}(x)+c_1\delta.
\end{split}
\end{equation}
Notice that by Proposition \ref{prop_U}, for every $x\ge a$, $U_a(x)$ is continuous as a function of $a$. Letting $\delta\to 0$ in \eqref{eq:asfeqn} we get
\begin{align}\notag
&\limsup\;\E^{Q^n}\left[\sup_{0\le t\le T} \int_0^{\tau_a[\tilde W^n]\wedge t}[l(\tilde W^n(s))-\lambda-\frac{1}{2}(\dot\psi^n(s))^2]ds+\rho 1_{\{{\tau_a[\tilde W^n]\le t\}}}\right]\\\notag
&\quad\le U_{a}(x)=U(x),
\end{align}
which is what we want to prove.

We now show that \eqref{refe3} holds. By \eqref{refe4a} and \eqref{refe4b}, we have for every $t\in[0,T]$ 
\begin{align}\notag
\tilde W^n(t)-\ph^n(t)=\int_0^t\Big[&-e(\tilde W^n(s))-e(\ph^n(s))+(\mu-r)(\pi^*(\tilde W^n(s))-\pi^*(\ph^n(s)))\\\notag
&+\sigma(\pi^*(\tilde W^n(s))-\pi^*(\ph^n(s)))\dot\psi^n(s)\Big]\;ds+\xi^n(t),
\end{align}
where
\begin{align}\notag
\xi^n(t):=\frac{1}{\sqrt{n}}\sigma\pi^*(\tilde W^n(t)dB^n(t),\quad t\in[0,T].
\end{align}
Since $e(\cdot)$ and $\pi^*(\cdot)$ are Lipschitz, there is a constant $c_2>0$ such that for every $t\in[0,T]$
\begin{align}\notag
|\tilde W^n(t)-\ph^n(t)|\le c_2\int_0^t(1+|\dot\psi^n(s)|)|\tilde W^n(s)-\ph^n(s)|ds+|\xi^n|_t
\end{align}
By Gronwall's inequality it follows that
\begin{align}\notag
|\tilde W^n(t)-\ph^n(t)|\le e^{c_2\int_0^t(1+|\dot\psi^n(s)|ds)}|\xi^n|_t.
\end{align}
Therefore,
\begin{align}\label{refe7}
|\tilde W^n-\ph^n|_{\tau_a[\tilde W^n]\wedge T}\le e^{c_2\int_0^T(1+|\dot\psi^n(s)|ds)}|\xi^n|_{\tau_a[\tilde W^n]\wedge T}.
\end{align}
For any $K>T$, consider 
\begin{align}\label{refe5}
B_n(K)=\left\{\omega: \int_0^T(1+|\dot\psi^n(s)|ds)\le K\right\}.
\end{align}
Clearly,
\begin{align}\notag
(A_n(\delta))^c=((A_n(\delta))^c\cap B_n(K))\cup((A_n(\delta))^c\cap (B_n(K))^c).
\end{align}
From \eqref{refe6}, \eqref{refe7}, and \eqref{refe5}, it follows that
\begin{align}\notag
(A_n(\delta))^c\cap B_n(K)\subset\left\{\omega: |\xi^n|_{\tau_a[\tilde W^n]\wedge T}\ge e^{-c_2 K}\delta\right\}.
\end{align}
By Doob's martingale inequality we have
\begin{align}\notag
\E^{Q^n}\Big[|\xi^n|_{\tau_a[\tilde W^n]\wedge T}^2\Big]\le4\E^{Q^n}\Big[\xi^n(\tau_a[\tilde W^n]\wedge T)^2\Big]\le \frac{c_2T}{n}.
\end{align}
Thus,
\begin{align}\notag
Q^n((A_n(\delta))^c\cap B_n(K))\le \frac{e^{2c_2K}}{\delta^2}\E^{Q^n}\Big[|\xi^n|_{\tau_a[\tilde W^n]\wedge T}^2\Big]\le \frac{c_2Te^{2c_2K}}{n\delta^2}.
\end{align}
On the other hand,
\begin{align}\notag
Q^n((A_n(\delta))^c\cap (B_n(K))^c)&\le Q^n((B_n(K))^c)\le Q^n\left(\int_0^T|\dot\psi^n(s)|^2ds>\frac{(K-T)^2}{T}\right)\\\notag
&\le\frac{T}{(K-T)^2}\E^{Q^n}\Big[\int_0^T|\dot\psi^n(s)|^2ds\Big],
\end{align}
where the second inequality follows from \eqref{refe5}. Due to \eqref{new7} and \eqref{cor1}, there is $N_1>0$ such that for every $n\ge N_1$, we have 
$$\E^{Q^n}\Big[\int_0^T|\dot\psi^n(s)|^2ds\Big]\le4\rho.$$
Hence, for every $n\ge N_1$, we have
\begin{align}\notag
Q^n((A_n(\delta))^c\cap (B_n(K))^c)&\le\frac{4\rho T}{(K-T)^2}.
\end{align}
Fix $\eps>0$. Let $K>T$ be such that 
$$\frac{4\rho T}{(K-T)^2}\le\frac{\eps}{2}.$$
Take $N=\max\{N_1,N_2\}$. then for every $n\ge N$, we have
$$Q^n((A_n(\delta))^c)<\eps.$$
This implies \eqref{refe3}.

\hfill$\Box$

\skp\noi
\textbf{Acknowledgement:} We thank the two anonymous referees, the AE and Huy\^{e}n Pham for insightful comments, which helped us improve our paper. We are also grateful to Virginia Young for many discussions that we had on the subject. This research is supported in part by the National Science Foundation through the DMS-1613170 grant.


%
%
%
\bibliographystyle{plain}

\bibliography{bc}

\end{document}